
\font\tenbb=msbm10     
\font\sevenbb=msbm7
\newfam\bbfam \def\bb{\fam\bbfam\tenbb}  
\textfont\bbfam=\tenbb
\scriptfont\bbfam=\sevenbb

\font\tenam=msam10     
\font\sevenam=msam7
\newfam\amfam   
\textfont\amfam=\tenam
\scriptfont\amfam=\sevenam

\mathchardef\eop="1903  

\font\bn=cmb10

\let\titlefont=\twelverm
\let\subheadfont=\bf

\let\em=\bf


\baselineskip=14pt \parskip=1.5pt \parindent=1.5em
\smallskipamount=4pt plus2pt minus1pt 
\medskipamount=7pt plus3pt minus2pt 
\bigskipamount=15pt plus5pt minus4pt 


\def\subhead#1{\goodbreak \removelastskip\penalty55\bigskip%
{\subheadfont #1}\par%
\nobreak\bigskip}


\def\intz{{\bb Z}}     
\def\GI{{\bb Z}[i\/]}    
\def\gcd{{\rm gcd}}    
\let\by\times


\long\def\profess#1#2\endprofess {\ifdim\lastskip<\medskipamount%
  \removelastskip\medskip\penalty-500\fi
  \noindent{\bn#1.\ }{\sl#2\par}%
  \ifdim\lastskip<\smallskipamount\removelastskip\penalty-55\smallskip\fi}

\long\def\rmprofess#1#2\endrmprofess {\ifdim\lastskip<\medskipamount%
  \removelastskip\medskip\penalty-350\fi
  \noindent{\bn#1.\ }{\rm#2\par}%
  \ifdim\lastskip<\smallskipamount\removelastskip\penalty-100\smallskip\fi}

\def\proof#1{\noindent {\it Proof#1.\enspace}}
\def\endproof{\ $\eop$\par%
   \ifdim\lastskip<\medskipamount
   \removelastskip\penalty-50\medskip\fi}

\def\({{\rm (}} \def\){{\rm )}}


\long\def\extrababble#1{\relax}
\def\comment#1{\relax}

\def\Descartessectionnum{4}
\def\uvquadruplepropnum{1}
\def\uvquadrupleprop{Proposition \uvquadruplepropnum%
\comment{ (uvquadrupleprop)}}
\def\quadruplethmnum{1}
\def\quadruplethm{Theorem \quadruplethmnum\comment{ (quadruplethm)}}
\def\uvsextuplethmnum{2}
\def\uvsextuplethm{Proposition \uvsextuplethmnum\comment{ (uvsextuplethm)}}
\def\sextuplethmnum{2}
\def\sextuplethm{Theorem \sextuplethmnum\comment{ (sextuplethm)}}
\def\uvquintuplethmnum{3}
\def\uvquintuplethm{Proposition \uvquintuplethmnum\comment{ (uvquintuplethm)}}
\def\quintuplethmnum{3}
\def\quintuplethm{Theorem \quintuplethmnum\comment{ (quintuplethm)}}
\def\descartesparamnum{4}
\def\descartesparamthm{Theorem \descartesparamnum\comment{ (descartesparamthm)}}

\def\ntupledefnum{1}
\def\uvntupledefnum{2}
\def\reversesubsnum{3}
\def\indianparanum{4}
\def\Cquadrnum{5}
\def\uvquadrupleeqnum{6}
\def\variantquadrsolutionnum{7}
\def\GLactionnum{8}
\def\quadruplematrixnum{9}
\def\sextupleeqnum{10}
\def\uvsextupleeqnum{11}
\def\polysextnum{12}
\def\Bnum{(13)}
\def\setonevarzeroinsextnum{14}
\def\DCnum{15}
\def\DCuvbijectionnum{(16)}

\line{\ \hfill\ }
\vskip-2em
\line{\tenrm To appear in J.~Pure Appl.~Algebra.\hfill}
\vskip2em

\centerline{\titlefont
Polynomial parametrization of Pythagorean 
quadruples, quintuples and sextuples.
}
\bigskip

\centerline{
Sophie Frisch (TU Graz) and Leonid Vaserstein (PSU)
}
\bigskip\bigskip

Abstract.
For $n=4$ or $6$, the Pythagorean $n$-tuples admit a parametrization
by a single $n$-tuple of polynomials with integer coefficients (which
is impossible for $n=3$). For $n=5$, there is an integer-valued
polynomial Pythagorean $5$-tuple which parametrizes Pythagorean
$5$-tuples (similar to the case $n=3$). Pythagorean quadruples are
closely related to (integer) Descartes quadruples, which we also
parametrize by a Descartes quadruple of polynomials with integer
coefficients.
\bigskip

\subhead{Introduction}

An Pythagorean triple is a triple of integers $(x_1,x_2,x_3)$ 
satisfying $x_1^2+ x_2^2=x_3^2$.
More generally, for any integer $n \ge 3$,  and any commutative ring $A$,
a Pythagorean $n$-tuple over $A$ is an $n$-tuple 
$(x_1,\ldots, x_n) \in A^n$ such that
$$ 
x_1^2 + \cdots + x_{n-1}^2 = x_n^2. \eqno(\ntupledefnum)
$$
Whenever $A$ is not specified, we will understand $A=\intz$. Likewise,
``polynomial Pythagorean $n$-tuple'' means a Pythagorean $n$-tuple 
over a ring of polynomials in finitely many indeterminates with
coefficients in $\intz$.

Instead of studying (\ntupledefnum) directly, it is often convenient to 
substitute $u=x_n+x_{n-1}$, $v=x_n-x_{n-1}$, and to consider the equation
$$
x_1^2 + \cdots + x_{n-2}^2 = uv. \eqno(\uvntupledefnum)
$$

Over any ring $A$ in which $2$ is not a zero-divisor, this substitution
and the reverse substitution 
$$x_{n-1}= (u-v)/2,\qquad x_n=(u+v)/2, \eqno(\reversesubsnum)$$
establish
a bijection between solutions $(x_1,\ldots, x_n)\in A^n$
of (\ntupledefnum) and solutions $(x_1,\ldots,x_{n-2},u,v)\in A^n$ with
$u-v\in 2A$  of (\uvntupledefnum).

We recall the existing polynomial parametrizations of integer Pythagorean
triples.
\extrababble{
Euclid in Book X, lemma 1 (circa 300 BC) gives a formula that produces
an infinite number of integer Pythagorean triples:
$$
(x_1,x_2,x_3)= (2y_1y_2, y_1^2-y_2^2, y_1^2+y_2^2). 
$$
}
It is well known that up to permutation of $x_1$ and $x_2$, every
Pythagorean triple has the form
$$ 
(x_1,x_2,x_3)= y_0(2y_1y_2, y_1^2-y_2^2, y_1^2+y_2^2) \eqno{(\indianparanum)}
$$
with  $y_i \in \intz$. In other words, the set of integer Pythagorean
triples is the union of $f_1(\intz^3)$ and  $f_2(\intz^3)$ where
$$
f_1(y_0, y_1, y_2) = y_0(2y_1y_2, y_1^2-y_2^2, y_1^2+y_2^2)  
$$
and
$$
f_2(y_0, y_1, y_2)= y_0(y_1^2-y_2^2, 2y_1y_2, y_1^2+y_2^2) 
$$
are two Pythagorean triples over the polynomial ring $\intz[y_0,y_1,y_2].$ 
We say that all Pythagorean triples are covered by two polynomial
Pythagorean triples (in 3 parameters each).

It is easy to see that the intersection of $f_1(\intz^3)$ and $f_2(\intz^3)$
contains only the zero triple. We know that it is not possible to cover all
Pythagorean triples by any one Pythagorean triple over
$\intz[y_1,\ldots, y_m]$ for any $m$ [FV]. It is, however, possible,
to cover all integer Pythagorean triples by a single Pythagorean triple 
over the ring of integer-valued polynomials in 4 indeterminates [FV].
An integer-valued polynomial is a polynomial with rational coefficients
which takes integer values whenever the variables take integer values.

\extrababble{
Since $f_1(y_0c^2,y_1,y_2) = f_1(y_0,y_1c,y_2c),$
different triples of parameters may give the same Pythagorean triple. 
The value for $y_0$ in the representation of a triple
$(x_1,x_2,x_3)\ne (0,0,0)$ 
could be made unique by the following convention:
$y_0 = \gcd(x_1,x_2,x_3){\rm sign}(x_3).$ 
Another reasonable choice for $y_0$ would be the square-free part
of $\gcd(x_1,x_2,x_3)$ times the sign of $x_3.$
Under either of these conventions for $y_0$, $(y_1,y_2)$ is unique up
to sign.
}

The primitive Pythagorean triples $(x_1,x_2,x_3)$  with positive $x_3,$
are, up to switching $x_1$ and $x_2,$  given by (\indianparanum) with
primitive $(y_1,y_2) \in \intz^2$ such that $y_1+y_2$ is odd.
The set of such pairs $(y_1,y_2) $ admits a polynomial parametrization [V].
Thus, all primitive Pythagorean triples can be covered by 4 polynomial
triples (in 95 parameters each, see [V], Example 14).

All positive Pythagorean triples are, up to switching of $x_1$ and $x_2$,
given by (\indianparanum) with integers $y_1>y_2>0$, $y_0 > 0.$
The set of such pairs admits a polynomial parametrization using the
fact that every positive integer can be written as a sum of 4 squares
plus $1$. Thus, the positive  Pythagorean triples
can be covered by 2 polynomial Pythagorean triples in 12 parameters.

It is unknown whether the set of positive primitive Pythagorean triples
can be parametrized by a finite set of polynomial Pythagorean triples.

\subhead{1. Quadruples}

After the short discussion of Pythagorean triples in the introduction,
we now address the case $n = 4$, in other words, Pythagorean quadruples.

Pythagorean quadruples were described by Carmichael [C], Chpt.~II, \S 10,
as follows:
up to permutation of $x_1,x_2,x_3$, every Pythagorean quadruple has the form 
$$
\eqalignno{
&(x_1,x_2,x_3,x_4)= f(y_0,y_1,y_2,y_3,y_4)=&\cr
&y_0(  2y_1y_3 + 2y_2y_4,\;  2y_1y_4 - 2y_2y_3,\;
y_1^2 + y_2^2 - y_3^2 - y_4^2,\;  y_1^2 + y_2^2 + y_3^2 + y_4^2 )&
(\Cquadrnum)\cr
}
$$
with integer values for the parameters $y_0,\ldots,y_4.$
Thus, all Pythagorean quadruples are covered by 6 polynomial
Pythagorean quadruples (in 5 parameters).
Considering the position of the odd entry, it is easy 
to see that at least 3 permutations of $x_1,x_2,x_3$ are needed.
If one examines Carmichael's proof, one sees that three polynomial
quadruples suffice, namely $(x_1,x_2,x_3,x_4)$,
$(x_1,x_3,x_2,x_4)$ and $(x_3,x_2,x_1,x_4)$.

We now show that a single polynomial Pythagorean quadruple covers all
Pythagorean quadruples. Our proof does not make use of Carmichael's
result (but rather provides a shorter proof of Carmichael's result as a 
byproduct). Nor do we use unique factorization in the ring of Gaussian
integers $\intz[i]$ (which could be used to give an alternative proof).

\rmprofess{Definition} 
An $n$-tuple $w= (w_1,\ldots, w_n)\in  A^n$ is called {\bf unimodular}
if $w_1A+\cdots + w_nA = A.$ In the case when $A = \intz$ this means
that $\gcd(w_1,\ldots, w_n) = 1,$ i.e., $w$ is primitive.
\endrmprofess

\profess{\uvquadrupleprop}
The integer solutions of
$$ x_1^2 + x_2^2 = uv  \eqno{(\uvquadrupleeqnum)} $$
are parametrized by the polynomial quadruple
$$
(x_1,x_2,u,v)=
y_0(y_1y_3+y_2y_4,\; y_1y_4-y_2y_3,\; y_1^2+y_2^2,\; y_3^2+y_4^2)
\eqno(\variantquadrsolutionnum)
$$
as the parameters vary through the integers. Also, $y_0$ can be
restriced to odd integers.
\endprofess

\proof{}
We represent the integer solutions of (\uvquadrupleeqnum) 
as Hermitian matrices
$$
w=\pmatrix{u &   x_1+i x_2  \cr x_1- i x_2  & v}
=\pmatrix{u &   x  \cr \bar x  & v}
$$
of determinant $0$ over the Gaussian integers $\GI$.

The group $GL(2, \GI)$ acts on the Hermitian matrices as follows
$$
w \to  g^*wg \eqno{(\GLactionnum)}
$$
where $g \in GL(2,\GI)$ and * means transposition composed with
entry-wise action of complex conjugation.  In particular, for an
elementary matix $g=E_{12}(\lambda)$ with $\lambda=\lambda_1+\lambda_2 i$

$$
\pmatrix{1 & 0  \cr \bar\lambda  & 1}
\pmatrix{u & x  \cr \bar x  & v}
\pmatrix{1 & \lambda  \cr 0  & 1} =
\pmatrix{u & (x_1+\lambda_1 u) + (x_2+\lambda_2 u)i\cr
(x_1+\lambda_1 u) - (x_2+\lambda_2 u)i& 
v +(\lambda_1^2+\lambda_2^2)u +2(\lambda_1x_1+\lambda_2x_2)}
$$

Setting either $\lambda_1=0$ or $\lambda_2=0$,  we see that we can 
add an arbitrary integer multiple of $u$ to $x_2$, leaving $u$ and
$x_1$ unchanged, and we can add an arbitrary multiple of $u$ to
$x_1$, leaving $u$ and $x_2$ unchanged. 

Given any solution with $u \ne 0$, we can, by applying elementary
matrices $g$ in $GL(2, \GI)$, make  $|x_1|, |x_2| \le |u|/2,$  and
hence $ \left|v\right| \le \left|u\right|/2.$   
Using the nontrivial permutation matrix in $GL(2, \GI)$, we can switch
$u$ and $v$.

Therefore, by an argument of descent, the orbit under $GL(2, \GI)$ of
any unimodular solution $w=(x_1,x_2,u,v)$ to (\uvquadrupleeqnum) contains
a solution with $v=0$ (and hence $x_1=x_2=0$) and $u = 1$ or $-1.$
So we get
$$
w = \pmatrix{u &   x_1+ x_2 i \cr x_1- x_2 i& v}
=  g^*\pmatrix{c & 0 \cr 0 & 0} g
=  \pmatrix{\bar a \cr \bar b} c (a,b)  \eqno(\quadruplematrixnum)
$$
with $c = \pm 1$, where $(a,b)$ is the first row of the matrix
$g \in GL(2,\GI).$

So every integer solution of (\uvquadrupleeqnum) has the form
$$
w=
\pmatrix{u &   x_1+ x_2 i  \cr
x_1- x_2 i  & v}
=  \pmatrix{\bar a \cr \bar b} c (a,b)  
$$
with $a,b\in\GI$, $c \in \intz$. Conversely, every expression of
this form is a solution to (\uvquadrupleeqnum) -- it is not necessary
to restrict $(a,b)$ to be primitive or $c$ to be $\pm 1$ or
(sum of 2 squares)-free.  

Writing $c=y_0$, $a=y_1+y_2i$ and $b=y_3+y_4i$,
with indeterminates $y_k$, we obtain a polynomial solution
$$
(u,v,x_1,x_2)=y_0(y_1^2+y_2^2,\; y_3^2+y_4^2,\; y_1y_3+y_2y_4,\; y_1y_4-y_2y_3)
$$
to (\uvquadrupleeqnum) which covers all integer solutions. 
If we replace $(a,b)$ above by $(1+i)(a,b)$, the solution is multiplied by 
$2$. We can, therefore, restrict $y_0$ to odd integers.
\endproof

\profess{\quadruplethm}
Let
$$
f(y_0,y_1,y_2,y_3,y_4)=
y_0(  2y_1y_3 + 2y_2y_4,\;  2y_1y_4 - 2y_2y_3,\;
y_1^2 + y_2^2 - y_3^2 - y_4^2,\;  y_1^2 + y_2^2 + y_3^2 + y_4^2 )
$$
The polynomial Pythagorean quadruple
$$
g=f(y_0/2,\; y_1,\; y_2,\; y_3,\; y_1+y_2+y_3+2z)\in \intz[y_0,y_1,y_2,y_3,z]
$$
in 5 parameters covers all Pythagorean quadruples, i.e., the range of
the function $g\colon \intz^5\rightarrow \intz^4$ consists of all
Pythagorean quadruples.
\endprofess

\proof{}
Applying (\reversesubsnum) to \uvquadrupleprop, we see
that every Pythagorean quadruple has the form
$$
(x_1,x_2,x_3,x_4) =  
y_0(  y_1y_3+y_2y_4,\; y_1y_4-y_2y_3,\; 
(y_1^2+y_2^2- y_3^2-y_4^2)/2,\; (y_1^2+y_2^2+ y_3^2+y_4^2)/2 )
$$
where $y_i \in \intz$ and $y_0(y_1+y_2+y_3+y_4)$ is even.
Since $y_0$ can be chosen odd, we may assume that $y_1+y_2+y_3+y_4$ is even.
Writing  $y_4=y_1+y_2+y_3+2z,$ we prove \quadruplethm.
\endproof

To get Carmichael's result, note that $x_4$ is odd for any primitive
Pythagorean quadruple $(x_1,x_2,x_3,x_4)$ and that exactly one of
$x_1,x_2,x_3$ is also odd.
So we can make $x_3+x_4$ even by switching, if necessary, $x_3$ with 
$x_1$ or $x_2.$
Then $\gcd(x_1,x_2,u,v) = 2$ for the corresponding solution
$(x_1,x_2,u,v)$ of (\uvquadrupleeqnum).
Going back from (\variantquadrsolutionnum) with $y_0=2$ to the
Pythagorean quadruple, we obtain Carmichael's formulas.

Notice that these formulas with $y_0 =1$ and primitive
$(y_1,y_2,y_3,y_4)$ do not necessary give primitive solutions. 
Our proof shows that the necessary and sufficient condition 
for primitivity is the primitivity of 
$(a,b) =(y_1+y_2i, y_3+y_4i)$. 
The set of primitive pairs of Gaussian integers admits a polynomial
parametrization by methods of [V], but this is beyond the scope of
the present paper.

\smallskip


\subhead{2. Sextuples}

We discuss Pythagorean sextuples before quintuples because we will
use sextuples in the proof of the parametrization of Pythagorean
quintuples by a single quintuple of integer-valued polynomials in
the next section.

Dickson [Di], Section 106, attempted to describe all Pythagorean  
sextuples, i.e., all integer solutions to 
$$
x_1^2+\ldots +x_5^2=x_6^2  \eqno{(\sextupleeqnum)}
$$
He observed that every integer solution of (\sextupleeqnum)
gives rise to an integer solution of
$$
x_1^2+\ldots +x_4^2=uv, \eqno{(\uvsextupleeqnum)}
$$
He solves the equation (\uvsextupleeqnum), in a lengthy proof of some
6 pages, but then fails to address the question of the reverse
substition: how to return to Pythagorean sextuples from integer
solutions of (\uvsextupleeqnum).

We will also start by parametrizing the integer solutions of
(\uvsextupleeqnum), giving a short proof using quaternions.

\rmprofess{Definition}
The algebra of {\em Lipschitz quaternions} is the $\intz$-algebra $L$
generated by two symbols $i$ and $j$ subject to the defining
relations $i^2= -1$, $j^2 = -1$, and $ji=-ij$. We set $k=ij$. 
\endrmprofess

We recall a few facts about the algebra of Lipschitz quaternions.
$L$ is a free $\intz$-module with basis $1, i, j, k$ and a free
$\GI$-module with basis $1,j$. $L$ can be represented as an
algebra of $4\times 4$ integer matrices or as an algebra of
$2\times 2$ matrices over $\GI$ by identifying $w=a+bi+cj+dk$
with
$$
M_4(w)=\pmatrix{\phantom{-}a&\phantom{-}b&\phantom{-}c&\phantom{-}d\cr
-b&\phantom{-}a&-d&\phantom{-}c\cr
-c&\phantom{-}d&\phantom{-}a&-b\cr 
-d&-c&\phantom{-}b&\phantom{-}a\cr}
\quad\hbox{\rm\ or\ }\quad
M_2(w)=\pmatrix{\phantom{-}a+bi&c+di\cr -c+di&a-bi\cr}
$$
respectively. 

An involution on $L$ is given by the $\intz$-algebra anti-isomorphism
$$a+bi+cj+dk\mapsto (a+bi+cj+dk)^*=a-bi-cj-dk.$$ In the $4\times 4$
integer matrix representation this corresponds to transposition;
and in the $2\times 2$ Gaussian integer matrix representation, to
transposition followed by complex conjugation. 

\rmprofess{Definition}
The {\em norm} of $w=a+bi+cj+dk \in L$ is defined as
$$(a^2+b^2+c^2+d^2)^2=\det( M_4(w)) = (\det M_2(w))^2$$
and the {\em reduced norm} as
$$a^2+b^2+c^2+d^2=w^*w=\det M_2(w).$$ 

For a $2\by 2$ matrix $M$ over $L$ we define the norm of $M$ as
the determinant of the $8\by 8$ integer matrix obtained by
replacing each matrix entry $w$ by $M_4(w)$, and the reduced norm
as the determinant of the $4\by 4$ matrix over $\GI$ obtained by
replacing each matrix entry $w$ by $M_2(w)$. 
\endrmprofess

\rmprofess{Remark} 
If $w$ is a Hermitian $2\by 2$ matrix over $L$, its entries
commute and we can calculate the determinant in a na\"\i ve way,
as 
$$\det \pmatrix{u & x_1+ x_2 i + x_3 j  +x_4 k \cr
x_1- x_2 i - x_3 j  -x_4 k  & v} = uv - x_1^2- x_2^2-x_3^2-x_4^2.
$$
The reduced norm of $w$ as defined above is the square
of this determinant.
\endrmprofess

\profess{\uvsextuplethm}
A parameterization of all integer solutions
$(x_1,x_2,x_3,x_4, u, v)$ of 
$$x_1^2+\ldots +x_4^2=uv$$
in $9$ parameters is given by
$$\eqalignno{
x_1&=y_0(y_1y_5+y_2y_6+y_3y_7+y_4y_8)\cr
x_2&=y_0(-y_1y_6+y_2y_5+y_3y_8-y_4y_7)\cr
x_3&=y_0(-y_1y_7-y_2y_8+y_3y_5+y_4y_6)\cr
x_4&=y_0(-y_1y_8+y_2y_7-y_3y_6+ y_4y_5)\cr
u&=y_0(y_1^2+y_2^2+y_3^2+y_4^2)\cr
v&=y_0(y_5^2+y_6^2+y_7^2+y_8^2),\cr
}$$
as the parameters $y_0,\ldots, y_8$ range through the integers. 
Here $y_0$ may be restriced to $\pm 1$. 
\endprofess

\proof{}
We identify integer solutions $w=(x_1,x_2,x_3,x_4,u,v)$ of
$x_1^2+\ldots +x_4^2=uv$ 
with $2\times 2$ Hermitian matrices over the algebra of 
Lipschitz quaternions $L$
$$
w=
\pmatrix{u & x_1+ x_2 i + x_3 j  +x_4 k \cr 
x_1- x_2 i - x_3 j  -x_4 k  & v}
$$
of reduced norm $0$.
 
The group $GL(2, L)$ acts on the set of $2\times 2$ Hermitian matrices
over $L$ of reduced norm $0$ by $(g,w)\mapsto g^*wg$, for $g\in GL(2, L)$,
where $g^*$ results from $g$ by application of the involution $*$ to
each entry, followed by transposition.

Given any unimodular solution of (\uvsextupleeqnum) with $u \ne 0$,
using an elementary matrix in $GL(2, L)$, we can make
$|x_1|, |x_2| , |x_3|, |x_4| \le |u|/2$, and hence $|v| \le |u|$.   
The inequality is strict unless $|x_m| = |v|/2$ for $m = 1,2,3,4$,
in which case $|v| =2$ by unimodularity.
In this last case, using an elementary matrix, we can arrange
$x_m =1$ for $m=1,2,3,4.$

Using the nontrivial permutation matrix in $GL(2, L)$, we can switch
$u$ and $v$.  Therefore, by induction on $|u|$, the orbit of any 
unimodular solution  $w=(x_1,x_2,x_3,x_4,u,v)$ to (\uvsextupleeqnum)
contains a solution with either $|u| = 1$ and $x_m = 0$
for $m = 1,2,3,4$  or $|u| = 2$ and  $x_m = 1$ for  $m = 1,2,3,4$.  

So we get that either
$$
w =  g^*\pmatrix{c & 0 \cr 0 & 0} g =   \pmatrix{a^* \cr b^*} c (a,b)  
$$
with $c = \pm 1$ where  $(a,b)$ is the first row of the matrix 
$g \in GL(2,L)$ or
$$
w= 
\pm g^* \pmatrix{2 & 1+ i +  j  +  k  \cr
1- i -  j  -  k  & 2} g 
= \pm \pmatrix{a^* \cr b^*} (a,b)  
$$
where $(a,b)=\pm (1- i, 1+ j )g$ with $g \in GL(2,L)$,
because
$$
\pmatrix{2 & 1+ i +  j  +  k  \cr 
1- i -  j  -  k  & 2}
= (1- i, 1+ j )^* (1- i, 1+ j ).
$$

So every integer solution of (\uvsextupleeqnum) has the form
$$
w=
\pmatrix{u & x_1+ x_2 i + x_3 j  +x_4 k \cr 
x_1- x_2 i - x_3 j  -x_4 k   & v}
=  \pmatrix{a^* \cr b^*} c (a,b)  
$$
with $a,b \in  L$ and $c \in\intz.$
Here we need not restrict $(a,b)$ to be primitive.

Writing 
$c=y_0, a = y_1+y_2  i+ y_3 j + y_4 k $
and $b = y_5+y_6  i+y_7 j+y_8 k $
we obtain the desired parameterization of all solutions 
$w=(x_1,x_2,x_3,x_4,u,v)$ of (\uvsextupleeqnum).

If we replace $(a,b)$ above by $d(a,b)$ with $d \in L$, the solution
$w$ is multiplied by $d^*d$, which is equivalent to replacing $y_0$ 
by $y_0d^*d$. Since every nonnegative integer is of the form $d^*d$
(sum of 4 squares) we can restrict $y_0$ to be $\pm 1.$
\endproof

Returning to Pythagorean $n$-tuples,
the following polynomial Pythagorean sextuple is known:
$$(x_1,\ldots,x_6)=h(y_0,\ldots,y_8)\in \intz[y_0,\ldots,y_8]^6 
\eqno(\polysextnum)$$
with
$$\eqalignno{
x_1 &= 2y_0 (y_1y_5+y_2y_6+y_3y_7+y_4y_8) \cr
x_2 &= 2y_0 (-y_1y_6+y_2y_5+y_3y_8-y_4y_7) \cr
x_3 &= 2y_0 (-y_1y_7-y_2y_8+y_3y_5+y_4y_6) \cr
x_4 &= 2y_0 (-y_1y_8+y_2y_7-y_3y_6+ y_4y_4) \cr
x_5 &= y_0 (y_1^2+y_2^2+y_3^2+y_4^2-y_5^2-y_6^2-y_7^2-y_8^2) \cr
x_6 &= y_0 (y_1^2+y_2^2+y_3^2+y_4^2+y_5^2+y_6^2+y_7^2+y_8^2) \cr
}
$$
There are, however, integer Pythagorean sextuples such as
$(1,1,1,1,0,2)$ that do not arise from the above polynomial sextuple
with integer parameters $y_i$. We now give a parametrization of
all integer Pythagorean sextuples by a single polynomial Pythagorean
sextuple in 9 parameters, or, by restricing the parameter $y_0$ to
$\pm 1$, a parametrization by two integer Pythagorean sextuples in
8 parameters each.

\profess{\sextuplethm}
Let $h=h(y_0,\ldots,y_8) \in  \intz[y_0,\ldots, y_8]^6$ be the 
polynomial Pythagorean sextuple \(\polysextnum\) above.
Then the polynomial Pythagorean sextuple in $\intz[y_0,\ldots, y_7,z]^6$
$$
g(y_0,y_1,y_2,y_3,y_4,y_5,y_6,y_7,z) =
h(y_0/2,y_1,y_2,y_3,y_4,y_5,y_6,y_7,y_1+y_2+y_3+y_4+y_5+y_6+y_7+2z) 
$$
in 9 parameters covers all Pythagorean sextuples, i.e., the range
of the function $g\colon \intz^9\rightarrow \intz^6$ is precisely
the set of all Pythagorean sextuples. Also, the parameter $y_0$
can be restricted to $\pm 1$.
\endprofess

\proof{}
We obtain all Pythagorean sextuples from all solutions of 
(\uvsextupleeqnum) with $u-v$ even by by (\reversesubsnum).
Since we we can take $y_0$ odd (even $\pm 1$) in \uvsextuplethm, we may
assume that $y_1+\cdots + y_8$ is even. 
Writing  $y_8=y_1+\cdots + y_7+2z,$
we obtain a Pythagorean sextuple over $\intz[y_0,\ldots,y_7,z]$
in 9 parameters which parametrizes all Pythagorean sextuples:

$$\eqalignno{
x_1&=y_0(y_1y_5+y_2y_6+y_3y_7+y_4(y_1+\cdots + y_7+2z))&\Bnum\cr
x_2&=y_0(-y_1y_6+y_2y_5+y_3(y_1+\cdots + y_7+2z)-y_4y_7)\cr
x_3&=y_0(-y_1y_7-y_2(y_1+\cdots + y_7+2z)+y_3y_5+y_4y_6)\cr
x_4&=y_0(-y_1(y_1+\cdots + y_7+2z)+y_2y_7-y_3y_6+ y_4y_5)\cr
x_5&=y_0(y_1^2+y_2^2+y_3^2+y_4^2- y_5^2-y_6^2-y_7^2
-(y_1+\cdots + y_7+2z)^2)/2\cr
x_6&=y_0(y_1^2+y_2^2+y_3^2+y_4^2+  y_5^2+y_6^2+y_7^2
+(y_1+\cdots + y_7+2z)^2)/2\cr
}$$
\endproof

\rmprofess{Remark}
The $\intz$-algebra $L$ is a subring of the ring $H$ of rational
quaternions (or Hamilton quaternions).  Adjoining to $L$ the element
$(1 + i + j + k)/2$, we obtain a larger subring $L'$ of $H$, called
the ring of {\em Hurwitz quaternions}. This ring $H$ has certain
unique factorization properties, which, however, we did not use.
They could be used to give alternative proofs for the results in 
this section.
\endrmprofess

\subhead{3. Quintuples}

We now consider the case $n = 5$ of Pythagorean quintuples.
We obtain \quintuplethm\  from \uvsextuplethm\ via \uvquintuplethm.

\profess{\uvquintuplethm}
A parametrization of all integer quintuples $(x_1, x_2, x_3, u, v)$ 
satisfying $$x_1^2+x_2^2+x_3^2=uv$$ by a quintuple of polynomials
with integer coefficients in the 12 parameters
$y_0,$ $z_0,$ $z_1,$ $z_2,$ $z_3,$ $z_4,$ $z_{12},$ $z_{13},$
$z_{14},$ $z_{23},$ $z_{24},$ $z_{34}$
is given by 
$$\eqalignno{
x_1&=y_0(y_1y_5+y_2y_6+y_3y_7+y_4y_8)\cr
x_2&=y_0(-y_1y_6+y_2y_5+y_3y_8-y_4y_7)\cr
x_3&=y_0(-y_1y_7-y_2y_8+y_3y_5+y_4y_6)\cr
u&=y_0(y_1^2+y_2^2+y_3^2+y_4^2)\cr
v&=y_0(y_5^2+y_6^2+y_7^2+y_8^2),\cr
}$$
with
$$\eqalignno{
y_1&=z_0z_1,\quad y_2=z_0z_2,\quad y_3=z_0z_3,\quad y_4=z_0z_4,\cr
y_5&=-z_{14}z_1-z_{24}z_2-z_{34}z_3\cr
y_6&=\phantom{-}z_{13}z_1+z_{23}z_2-z_{34}z_4\cr
y_7&=-z_{12}z_1+z_{23}z_3+z_{24}z_4\cr
y_8&=-z_{12}z_2-z_{13}z_3-z_{14}z_4.\cr
}$$
\endprofess

\proof{}
To solve 
$x_1^2+x_2^2+x_3^2=uv,$
we set $x_4=0$
in the general solution to
$x_1^2+x_2^2+x_3^2+x_4^2=uv$
obtained in \uvsextuplethm.
$$
-y_1y_8+y_2y_7-y_3y_6+ y_4y_5 = 0 \eqno(\setonevarzeroinsextnum)
$$ 
(The case $y_0=0$ only contributes the zero solution which 
we will not miss.)

The integer solutions of (\setonevarzeroinsextnum) can be
parametrized by 11 parameters as follows.

First we write  $(y_1,y_2,y_3,y_4) = z_0(z_1,z_2,z_3,z_4)$
with $z_i \in \intz$ and $(z_1,z_2,z_3,z_4)$ unimodular such that
$(z_1,z_2,z_3,z_4).(-y_8,y_7,-y_6,y_5) = 0$ (the case $y_1=y_2=y_3=y_4=0$
only contributes solutions $(0,0,0,0,v)$ which we will not miss).

By [VS], Remark after Lemma 9.6, we can write 
$$
\displaylines{
(-y_8,y_7,-y_6,y_5) =\hfill\cr
z_{12}(z_2,-z_1,0,0) +  z_{13}(z_3,0, -z_1,0) +  z_{14}(z_4,0,0, -z_1)
+  z_{23}(0,z_3,-z_2,0) +  z_{24}(0,z_4,0,-z_2) +  z_{34}(0,0,z_4,-z_3).
}
$$

This gives a parametrization of the integer solutions of 
(\setonevarzeroinsextnum) in the 11 parameters 
$z_0,$ $z_1,$ $z_2,$ $z_3,$ $z_4,$ $z_{12},$ $z_{13},$ 
$z_{14},$ $z_{23},$ $z_{24},$ $z_{34}.$

Therefore all integer solutions of 
$x_1^2+x_2^2+x_3^2=uv$
are parametrized by a polynomial solution with 12 parameters
including $y_0$.
\endproof

Another parametrization of $x_1^2+x_2^2+x_3^2=uv$ with 20 parameters
can be obtained using [V], Proposition 3.4 with $k=8$. We now parametrize
integer Pythagorean quintuples by a single Pythagorean quintuple
over the ring of integer-valued polynomials in 14 variables. This 
can be used to construct a parametrization by a finite number of
integer-coefficient polynomial Pythagorean quintuples [F]. Whether
it is possible to parametrize integer Pythagorean quintuples by a
single quintuple of integer-coefficient polynomials or not, we do
not know.

\profess{\quintuplethm}
A parametrization of all Pythagorean quintuples by a quintuple of
integer-valued polynomials in the 14 variables
$w_0, w_{12}, w_{13}, w_{14}, w_{23}, w_{24}, w_{34},
t_1, t_2, t_3, d_1, d_2, d_3, w_4$ 
is given by
$(f_1, f_2, f_3, f_5, f_6)$,
where
$$\eqalignno{
f_1 &= 2y_0 (y_1y_5+y_2y_6+y_3y_7+y_4y_8) \cr
f_2 &= 2y_0 (-y_1y_6+y_2y_5+y_3y_8-y_4y_7) \cr
f_3 &= 2y_0 (-y_1y_7-y_2y_8+y_3y_5+y_4y_6) \cr
f_5 &= y_0 (y_1^2+y_2^2+y_3^2+y_4^2-y_5^2-y_6^2-y_7^2-y_8^2)/2 \cr
f_6 &= y_0 (y_1^2+y_2^2+y_3^2+y_4^2+y_5^2+y_6^2+y_7^2+y_8^2)/2 \cr
} $$
and
$$\eqalignno{
y_1&=z_0z_1,\quad y_2=z_0z_2,\quad y_3=z_0z_3,\quad y_4=z_0z_4,\cr
y_5&=-z_{14}z_1-z_{24}z_2-z_{34}z_3\cr
y_6&=\phantom{-}z_{13}z_1+z_{23}z_2-z_{34}z_4\cr
y_7&=-z_{12}z_1+z_{23}z_3+z_{24}z_4\cr
y_8&=-z_{12}z_2-z_{13}z_3-z_{14}z_4.\cr
}$$
and further
\goodbreak
$$\eqalignno{
z_0=&
w_0+t_1 w_0+t_2 w_0-2 t_1 t_2 w_0+t_3 w_0-2 t_1 t_3 w_0-t_2 t_3 w_0
+2 t_1 t_2 t_ 3 w_0+t_1 w_{12}- t_1 t_2 w_{12}-\cr
&-t_1 t_3 w_{12}+t_2 t_3 w_{12}+t_2 w_{13}-t_1 t_2 w_{13}+t_3 w_{14}
-t_1 t_3 w_{14}+ t_1 w_{23} + t_2 w_{23}-2 t_1 t_2 w_{23}-\cr
&-t_1 t_3 w_{23} -t_2 t_3 w_{23} +2 t_1 t_2 t_3 w_{23} +t_1 w_{24}
-t_1 t_2 w_{24}+ t_3 w_{24}-2 t_1 t_3 w_{24}-t_2 t_3 w_{24}+\cr
&+2 t_1 t_2 t_3 w_{24} +t_2 w_{34}-t_1 t_2 w_{34}
+t_3 w_{34}- t_1 t_3 w_{34}-2 t_2 t_3 w_{34}+ 2 t_1 t_2 t_3 w_{34}\cr
z_1=& 2 d_1+t_1 t_2+t_3-2 t_1 t_2 t_3+w_4 \cr
z_2=&2 d_2+t_1-t_1 t_2+t_3-t_1 t_3-t_2 t_3+2 t_1 t_2 t_3+w_4 \cr
z_3=&2 d_3+t_2+t_3-t_1 t_3-2 t_2 t_3+2 t_1 t_2 t_3+w_4\cr
z_4=&w_4\cr
z_{12}=&w_{12}+t_1 t_2 w_{12}-t_1 t_2 t_3 w_{12}+t_1 t_2 w_{14}
-t_1 t_2 t_3 w_{14}+t_1 t_2 w_{23}-t_1 t_2 t_3 w_{23}+t_1 t_2 w_{34}
-t_1 t_2 t_3 w_{34}
\cr
z_{13}=&
w_{13}+t_1 t_3 w_{13}-t_1 t_2 t_3 w_{13}+t_1 t_3 w_{14}-t_1 t_2 t_3 w_{14}
+t_1 t_3 w_{23}-t_1 t_2 t_3 w_{23}+ t_1 t_3 w_{24}-t_1 t_2 t_3 w_{24}
\cr
z_{14}=& w_{14}\cr
z_{23}=& w_{23}\cr
z_{24}=&
t_1t_2 t_3 w_{12}+t_1t_2 t_3 w_{13}+w_{24}+t_1 t_2 t_3 w_{24}
+t_1 t_2 t_3 w_{34}
\cr
z_{34}=& w_{34}\cr
}$$
\endprofess

\proof{}
To go from the solutions of $x_1^2+x_2^2+x_3^2=uv$ parametrized in
\uvquintuplethm\ to the solutions of $x_1^2+x_2^2+x_3^2+x_4^2=x_5^2$
we use (\reversesubsnum), allowing only those $u,v$ with $u\pm v$ even.

In our case, we need to parametrize those $z_0,\ldots,z_{34}$ that make
$y_1+y_2+\ldots+y_8$ even, i.e., those $z_0,\ldots,z_{34}$ such that
$E= z_0z_1 + z_0z_2 + z_0z_3 + z_0z_4 +
-z_{14}z_1-z_{24}z_2-z_{34}z_3
-z_{13}z_1-z_{23}z_2+z_{34}z_4
-z_{12}z_1+z_{23}z_3+z_{24}z_4
-z_{12}z_2-z_{13}z_3-z_{14}z_4$
is even.

This is acheived by the following substitution, which, after
simplification, gives the parametrization in the statement of the theorem.
$$
\displaylines{
\bigl(z_0, z_1, z_2, z_3, z_4, 
z_{12}, z_{13}, z_{14}, z_{23}, z_{24}, z_{34}\bigr) =\cr
%
\bigl(w_0, w_4+2d_1, w_4+2d_2, w_4+2d_3, w_4,
w_{12}, w_{13}, w_{14}, w_{23}, w_{24}, w_{34}\bigr)(1-t_1)(1-t_2)(1-t_3) +\cr
%
\bigl(w_{14}+w_{24}+w_{34}+2w_0, w_4+2d_1+1, w_4+2d_2+1, w_4+2d_3+1, w_4,
w_{12}, w_{13}, w_{14}, w_{23}, w_{24}, w_{34}\bigr)(1-t_1)(1-t_2)t_3 +\cr
%
\bigl(w_{13}+w_{23}+w_{34}+2w_0, w_4+2d_1, w_4+2d_2, w_4+2d_3+1, w_4,
w_{12}, w_{13}, w_{14}, w_{23}, w_{24}, w_{34}\bigr)(1-t_1)t_2(1-t_3) +\cr
%
\bigl(w_{12}+w_{23}+w_{24}+2w_0, w_4+2d_1,w_4+2d_2+1, w_4+2d_3, w_4,
w_{12}, w_{13}, w_{14}, w_{23}, w_{24}, w_{34}\bigr)t_1(1-t_2)(1-t_3) +\cr
%
\bigl(w_{12}+w_{13}+w_{14}+2w_0, w_4+2d_1+1, w_4+2d_2, w_4+2d_3, w_4,
w_{12}, w_{13}, w_{14}, w_{23}, w_{24}, w_{34}\bigr)(1-t_1)t_2t_3 +\cr
%
\bigl(w_0, w_4+2d_1+1, w_4+2d_2+1, w_4+2d_3, w_4, 
w_{12}, w_{23}+w_{24}+w_{14}+2w_{13}, w_{14}, w_{23}, w_{24}, w_{34}\bigr)
t_1(1-t_2)t_3 +\cr
%
\bigl(w_0, w_4+2d_1+1, w_4+2d_2, w_4+2d_3+1, w_4,
w_{23}+w_{14}+w_{34}+2w_{12}, w_{13}, w_{14}, w_{23}, w_{24}, w_{34}\bigr)
t_1t_2(1-t_3) +\cr
%
\bigl(w_0, w_4+2d_1, w_4+2d_2+1, w_4+2d_3+1, w_4, 
w_{12}, w_{13}, w_{14}, w_{23}, w_{12}+w_{13}+w_{34}+2w_{24}, w_{34}\bigr)
t_1t_2t_3\cr
}
$$
\endproof

\subhead{\Descartessectionnum. Descartes quadruples}

In 1643 Descartes [De] described a relationship between the radii of four
mutually tangent circles (called a Descartes configuration), namely,
$$
2(b_1^2+b_2^2 + b_3^2 + b_4^2)=(b_1+b_2+b_3+b_4)^2 \eqno(\DCnum)
$$
where $b_i$ are the reciprocals of the radii. Others, including
Steiner, Beecroft, and Soddy [S], rediscovered the result.
We call an integer solution of (\DCnum) a Descartes quadruple.

Given one Descartes configuration, there is a geometric way to produce
plenty of them creating an Apollonian packing.
If the four intitial curvatures $b_i$ are integers, all curvatures in
the packing are integers.
There are several publications about integer Apollonian packings
[GLMWY1], [GLMWY2], [GLMWY3], [EL],[N], [LMW]. A bijection between
integer Pythagorean quadruples and integer Descartes quadruples
can be found in [GLMWY1], Lemma 2.1. 

In this section we parametrize all integer solutions of (\DCnum) by 
a single polynomial solution in 5 parameters, using a bijection between 
integer Descartes quadruples and integer solutions of (\uvquadrupleeqnum).

Given an integer solution $(x_1, x_2, u, v)$ of (\uvquadrupleeqnum),
$$
b_1 = u + v - 2 x_1 + x_2,\quad b_2 = u + x_2,\quad
b_3 = v + x_2,\quad b_4 = -x_2
\eqno{\DCuvbijectionnum}
$$
is an integer solution of (\DCnum). Conversely, we can invert this linear
transformation: given an integer solution $(b_1,b_2,b_3,b_4)$ of (\DCnum),
$b_1 + b_2 + b_3 + b_4$ is even and
$$
x_1= (-b_1 + b_2 + b_3 + b_4)/2,\quad x_2 = -b_4,\quad
u = b_2 + b_4,\quad v = b_3 + b_4
$$
is an integer solution of (\uvquadrupleeqnum).

\profess{\descartesparamthm}
A parametrization 
of all integer solutions 
of
$$
2(b_1^2+b_2^2 + b_3^2 + b_4^2)=(b_1+b_2+b_3+b_4)^2 \eqno(\DCnum)
$$
in 5 parameters is given by
$$
\eqalignno{
b_1 &= 
y_0(y_1^2+y_2^2 + y_3^2+y_4^2 - 2y_1y_3 -2y_2y_4 + y_1y_4-y_2y_3),\cr
b_2 &=
y_0(y_1^2+y_2^2 + y_1y_4-y_2y_3),\cr
b_3 &=
y_0(y_3^2+y_4^2 + y_1y_4-y_2y_3),\cr
b_4 &=
y_0(-y_1y_4-y_2y_3).\cr
}
$$
\endprofess

\proof{}
In the expression \DCuvbijectionnum\ of $b_1,b_2,b_3,b_4$ in terms of
a solution $x_1,x_2,u,v$ of (\uvquadrupleeqnum)
we have substituted the parametrization 
of all integer solutions of (\uvquadrupleeqnum) from
\uvquadrupleprop.
\endproof

\subhead{References}

[C]  R.D. Carmichael, Diophantine Analysis, New York: John Wiley \& Sons, 1915.
\smallskip

[CS]  John H. Conway and Derek A. Smith. 
On quaternions and octonions. Their geometry, arithmetic, and symmetry.
Peters, 2003.
\smallskip

[De]
R. Descartes, Oeuvres de Descartes, Correspondance IV, 
(C. Adam and P. Tannery, Eds.), Paris: Leopold Cert 1901.
\smallskip

[Di]  Dickson, Leonard Eugene, Algebras and their arithmetics.  
Dover Publications, Inc., New York 1960. 

[EL]
Eriksson, Nicholas; Lagarias, Jeffrey C.
Apollonian circle packings: number theory.
II. Spherical and hyperbolic packings.
Ramanujan J. 14 (2007), no. 3, 437--469. 
\smallskip

[F]
Sophie Frisch, Remarks on polynomial parametrization of sets 
of integer points,  Comm. Algebra  36  (2008),  no. 3, 1110--1114.
\smallskip

[FV]
Sophie Frisch and Leonid Vaserstein, 
Parametrization of {P}ythagorean triples
by a single triple of polynomials,
J.~Pure Appl.~Algebra 212 (2008) 271--274.
\smallskip

[GLMWY1]
R.L. Graham, J.C. Lagarias, C.L. Mallows, A.R. Wilks, C.R. Yan,
Apollonian circle packings: Number theory, J. Number Theory 100 (2003) 1--45. 
\smallskip

[GLMWY2]
Graham, Ronald L.; Lagarias, Jeffrey C.; Mallows, Colin L.; 
Wilks, Allan R.; Yan, Catherine H. 
Apollonian circle packings: geometry and group theory.
II. Super-Apollonian group and integral packings.
Discrete Comput. Geom. 35 (2006), no. 1, 1--36. 
\smallskip

[GLMWY3]
Graham, Ronald L.; Lagarias, Jeffrey C.; Mallows, Colin L.;
Wilks, Allan R.; Yan, Catherine H.
Apollonian circle packings: geometry and group theory. 
III. Higher dimensions.  
Discrete Comput. Geom. 35 (2006), no. 1, 37--72. 
\smallskip

[LMW]
J.C. Lagarias, C.L. Mallows, A.R. Wilks,
Beyond the Descartes circle theorem,
Amer. Math. Monthly 109 (4) (2002) 338--361. 
\smallskip

[N]
Northshield, S.
On integral Apollonian circle packings.  
J. Number Theory 119 (2006), no. 2, 171--193. 
\smallskip

[S] 
F. Soddy, The kiss precise, Nature, vol 137 no. 3477 (June 20 1936) p1021.

[S] Suslin, A. A.
The structure of the special linear group over rings of polynomials. (Russian)
Izv. Akad. Nauk SSSR Ser. Mat. 41 (1977), no. 2, 235--252, 477. 
\smallskip

[VS] L.N. Vaserstein and A.A. Suslin, Serre's problem on projective 
modules over polynomial rings and  algebraic $K$-theory,
Izv.Akad.Nauk, ser.mat. 40:5 (1976), 993-1054 = Math.USSR Izv. 10:5, 937-1001.
\smallskip

[V] L.N. Vaserstein,  Polynomial parametrization for the solutions
of Diophantine equations and arithmetic groups,
Annals of Math. 171 (2), 2010, 979--1009.

\bye